# Estimates for the difference between approximate and exact solutions to stochastic differential equations in the G-framework

Faiz Faizullah [a,b], Ilyas Khan [c,d,e], Mukhtar M. Salah [d,e] and Ziyad Ali Alhussain [e]

[a]Department of Mathematics, Swansea University, Swansea, UK; [b]Department of Basic Sciences and Humanities, College of Electrical and Mechanical Engineering, National University of Sciences and Technology (NUST) Islamabad Pakistan; [c]Faculty of Mathematics and Statistics, Ton Duc Thang University, Ho Chi Minh City, Vietnam; [d]Basic Engineering Sciences Department College of Engineering, Majmaah University, Majmaah, Saudi Arabia; [e]Department of Mathematics, College of Science Majmaah University, Majmaah, Saudi Arabia

**ABSTRACT**

This article investigates the Euler-Maruyama approximation procedure for stochastic differential equations in the framework of G-Browinian motion with non-linear growth and non-Lipschitz conditions. The results are derived by using the Burkholder-Davis-Gundy (in short BDG), Hölder's, Doobs martingale's and Gronwall's inequalities. Subject to non-linear growth condition, it is revealed that the Euler-Maruyama approximate solutions are bounded in $M_G^2([t_0, T]; \mathbb{R}^n)$. In view of non-linear growth and non-uniform Lipschitz conditions, we give estimates for the difference between the exact solution $Z(t)$ and approximate solutions $Z^q(t)$ of SDEs in the framework of G-Brownian motion.

## 1. Introduction

The stochastic differential equations (SDEs) theory has been used in several disciplines of sciences and engineering. In biological sciences, they are utilized to model the achievement of stochastic changes in reproduction on population processes [1,2]. In space, SDEs describe the transport of cosmic rays. They can be used to model the climate and weather. The percolation of fluid through absorbent structures and water catchment can be modelled by SDEs [3]. They are now very common in mechanical, computer, chemical and electrical engineering etc. By virtue of the growth and Lipschitz conditions, SDEs in the framework of G-Brownian motion were studied by Peng [4,5]. He derived the existence and uniqueness results in view of the contraction principle technique. With Picard's iteration scheme, the stated theory was developed by Gao [6]. By virtue of the Caratheodory approximation procedure, the existenc-uniqueness results were achieved by Faizullah [7]. The mentioned theory was extended to integral Lipschitz conditions by Bai and Lin [8]. Subject to the discontinuous coefficients, Faizullah derived that SDEs in the G-framework possess more than one solutions [9]. In the present article, we investigate the Euler-Maruyama approximation procedure for SDEs in the framework of G-Browinian motion with non-linear growth and non-Lipschitz conditions. Let $0 \leq t_0 \leq t \leq T < \infty$ and consider the following stochastic differential equation in the G-framework

$$dZ(t) = g(t, Z(t))\,dt + h(t, Z(t))\,d\langle W, W\rangle(t) + w(t, Z(t))\,dW(t), \quad (1)$$

on $t \in [t_0, T]$ with given initial condition $Z(t_0) = Z_0$. The quadratic variation process of G-Browniam motion $\{W(t) : t \geq 0\}$ is denoted by $\{\langle W, W\rangle(t) : t \geq 0\}$. For all $x \in \mathbb{R}^n$, the given coefficients $g(., x), h(., x)$ and $w(., x)$ belongs to space $M_G^2([t_0, T]; \mathbb{R}^n)$. The integral form of Equation (1) is given as the following

$$Z(t) = Z_0 + \int_{t_0}^{t} g(s, Z(s))\,ds + \int_{t_0}^{t} h(s, Z(s))\,d\langle W, W\rangle(s) + \int_{t_0}^{t} w(s, Z(s))\,dW(s), \quad t \in [0, T], \quad (2)$$

where $Z_0 \in \mathbb{R}^d$ is a given initial condition. All through the present article, we assume the following assumptions. Let $t \in [t_0, T]$. For every $u, v \in \mathbb{R}^n$

$$|g(t, u) - g(t, v)|^2 + |h(t, u) - h(t, v)|^2 + |w(t, u) - w(t, v)|^2 \leq \Upsilon(|u - v|^2), \quad (3)$$

where the function $\Upsilon(.) : \mathbb{R}^+ \to \mathbb{R}^+$ is non-decreasing and concave with $\Upsilon(0) = 0, \Upsilon(r) > 0$ for $r > 0$ and

$$\int_{0+} \frac{dr}{\Upsilon(r)} = \infty. \quad (4)$$

CONTACT Ilyas Khan  ilyaskhan@tdt.edu.vn

Since $\Upsilon$ is concave and $\Upsilon(0) = 0$, for all $r \geq 0$,

$$\Upsilon(r) \leq C + Dr, \tag{5}$$

where $C$ and $D$ are positive constants. For every $t \in [t_0, T]$ and $g(t, 0), h(t, 0), w(t, 0) \in L^2$,

$$|g(t, 0)|^2 + |h(t, 0)|^2 + |w(t, 0)|^2 \leq M, \tag{6}$$

where $M$ is a positive constant. Assumptions (3) and (6) are known as non-uniform Lipschitz and weakened linear growth conditions respectively. The current paper is organized in three more sections. Section 2 contains several basic results and concepts such as the definitions of G-expectation, G-Brownian motion, Ito's integral, Hölder's inequality, Doobs martingale's inequality and Gronwall's inequality etc. Section 3 presents the idea of Euler-Maruyama approximate solutions for SDEs in the G-framework. This section consists of an important result, which reveals that the Euler-Maruyama approximate solutions are bounded. In Section 4, we prove an important lemma, which is utilized in the main theorem. This section gives estimates for the difference between an exact solution $Z(t)$ and approximate solutions $Z^q(t)$ of SDEs in the framework of G-Brownian motion.

## 2. Preliminaries

Building on the previous ideas of G-Brownian motion theory, this section is devoted to the preliminary notions and results required for the subsequent sections of this article. For more details on the concepts of G-Brownian theory, readers are suggested to consult the papers [10–17]. Let $\Omega$ be a given fundamental non-empty set. Suppose $\mathcal{H}$ be a space of linear real functions defined on $\Omega$ satisfying that (i) $1 \in \mathcal{H}$ (ii) for every $d \geq 1, Z_1, Z_2, \ldots, Z_n \in \mathcal{H}$ and $\varphi \in C_{b.Lip}(\mathbb{R}^n)$ it holds $\varphi(Z_1, Z_2, \ldots, Z_n) \in \mathcal{H}$ i.e., with respect to Lipschitz bounded functions, $\mathcal{H}$ is stable. Then $(\Omega, \mathcal{H}, E)$ is a sub-expectation space, where $E$ is a sub-expectation defined as the following.

**Definition 2.1:** A functional $E : \mathcal{H} \to \mathbb{R}$ satisfying the below four features is known as a sub-expectation. Let $X, Y \in \mathcal{H}$, then

(i) Monotonicity: $E[Z] \geq E[Y]$ if $Z \geq Y$.
(ii) Constant preservation: $E[K] = K$, for all $K \in \mathbb{R}$.
(iii) Positive homogeneity: $E[\alpha Z] = \alpha E[Z]$, for all $\alpha \in \mathbb{R}^+$.
(iv) Sub-additivity: $E[Z] + E[Y] \geq E[Z + Y]$.

Moreover, let $\Omega$ be the space of all $\mathbb{R}^n$-valued continuous paths $(w_t)_{t \geq 0}$ starting from zero. Also, suppose that associated with the below distance, $\Omega$ is a metric space

$$\rho(w^1, w^2) = \sum_{i=1}^{\infty} \frac{1}{2^i} (\max_{t \in [0,k]} |w_t^1 - w_t^2| \wedge 1).$$

Fix $T \geq 0$ and set

$$L_{ip}^0(\Omega_T) = \{\phi(W_{t_1}, W_{t_2}, \ldots, W_{t_m}) : m \geq 1,$$
$$t_1, t_2, \ldots, t_m \in [0, T], \phi \in C_{b.Lip}(\mathbb{R}^{m \times n}))\},$$

where $W$ is the canonical process, $L_{ip}^0(\Omega_t) \subseteq L_{ip}^0(\Omega_T)$ for $t \leq T$ and $L_{ip}^0(\Omega) = \cup_{n=1}^{\infty} L_{ip}^0(\Omega_n)$. The completion of $L_{ip}^0(\Omega)$ under the Banach norm $E[|.|^p]^{1/p}$, $p \geq 1$ is denoted by $L_G^p(\Omega)$, where $L_G^p(\Omega_t) \subseteq L_G^p(\Omega_T) \subseteq L_G^p(\Omega)$ for $0 \leq t \leq T < \infty$. Generated by the canonical process $\{W(t)\}_{t \geq 0}$, the filtration is symbolized as $\mathcal{F}_t = \sigma\{W_s, 0 \leq s \leq t\}$, $\mathcal{F} = \{\mathcal{F}_t\}_{t \geq 0}$. Suppose $\pi_T = \{t_0, t_1, \ldots, t_N\}$, $0 \leq t_0 \leq t_1 \leq \cdots \leq t_N \leq \infty$ be a partition of $[0, T]$. Set $p \geq 1$, then $M_G^{p,0}(0, T)$ indicates a collection of the below type processes

$$\mu_t(w) = \sum_{i=0}^{N-1} \delta_i(w) I_{[t_i, t_{i+1}]}(t), \tag{7}$$

where $\delta_i \in L_G^p(\Omega_{t_i}), i = 0, 1, \ldots, N - 1$. Furthermore, the completion of $M_G^{p,0}(0, T)$ with the below given norm is indicated by $M_G^p(0, T), p \geq 1$

$$\|\mu\| = \left\{ \int_0^T E[|\mu_s|^p] \, ds \right\}^{1/p}.$$

**Definition 2.2:** A d-dimensional stochastic process $\{W(t)\}_{t \geq 0}$ satisfying the below properties is called a G-Brownian motion

(1) $W(0) = 0$.
(2) The increment $W_{t+m} - W_t$, for any $t, m \geq 0$, is G-normally distributed and independent from $W_{t_1}, W_{t_2}, \ldots \ldots W_{t_n}$, for $n \in N$ and $0 \leq t_1 \leq t_2 \leq, \ldots, \leq t_n \leq t$,

**Definition 2.3:** Let $\mu_t \in M_G^{2,0}(0, T)$ having the form (7). Then the G-quadratic variation process $\{\langle W \rangle_t\}_{t \geq 0}$ and G-Itô's integral $I(\mu)$ are respectively defined by

$$\langle W \rangle_t = W_t^2 - 2 \int_0^t W_u \, dW_u,$$

$$I(\mu) = \int_0^T \eta_u \, dW_u = \sum_{i=0}^{N-1} \delta_i(W_{t_{i+1}} - W_{t_i}).$$

The following two lemmas are borrowed from the book [18]. They are called as Hölder's and Gronwall's inequalities respectively.

**Lemma 2.4:** *Let $p, q > 1$, $1/p + 1/q = 1$ and $g, h \in L^2$. Then $gh \in L^1$ and*

$$\int_a^b g(t)h(t)\,dt \leq \left(\int_a^b |g(t)|^p\,dt\right)^{1/p} \left(\int_a^b |h(t)|^q\,dt\right)^{1/q}.$$

**Lemma 2.5:** *Let $g(t) \geq 0$ and $h(t) \geq 0$ be continuous real functions defined on $[a, b]$. If for all $t \in [a, b]$,*

$$h(t) \leq M + \int_a^b g(s)h(s)\,ds,$$

*where $M \geq 0$, then*

$$h(t) \leq M e^{\int_a^t g(s)\,ds},$$

*for all $t \in [a, b]$.*

For more details of the following (Burkholder-Davis-Gundy (BDG) inequalities) two lemmas, see [6].

**Lemma 2.6:** *Let $\eta \in M_G^p(0, T)$ then for any $p \geq 1$,*

$$\hat{\mathbb{E}}\left[\sup_{0 \leq v \leq T} \left|\int_0^v \eta(s)\,d\langle B\rangle(s)\right|^p\right] \leq k_1 T^{p-1} \hat{\mathbb{E}} \int_0^T |\eta(s)|^p\,ds,$$

*where $0 < k_1 < \infty$ is a positive constant depends only on $p$.*

**Lemma 2.7:** *Let $\eta \in M_G^2(0, T)$ then for any $p \geq 2$,*

$$\hat{\mathbb{E}}\left[\sup_{0 \leq v \leq T} \left|\int_0^v \eta(s)\,dB(s)\right|^p\right] \leq k_2 T^{p/2-1} \hat{\mathbb{E}} \int_0^T |\eta(s)|^p\,ds,$$

*where $0 < k_2 < \infty$ is a positive constant depends only on $p$.*

Just for simplicity, all through this article we take $k_1 = k_2 = 1$.

## 3. Euler-Maruyama approximate solutions

We now describe Euler-Maruyama approximation procedure for Equation (2). For any integer $q \geq 1$, we define a sequence $\{Z^q(t)\}$ on $[t_0, T]$ as follows. For $t \in [t_0, t_0 + n/q]$, $Z^q(t_0) = Z_0$ and for $t_0 + n/q < t \leq (t_0 + (n+1)/q) \wedge T, n = 1, 2, \ldots$,

$$Z^q(t) = Z^q\left(t_0 + \frac{n}{q}\right) + \int_{t_0+n/q}^t g\left(s, Z^q\left(t_0 + \frac{n}{q}\right)\right) ds$$
$$+ \int_{t_0+n/q}^t h\left(s, Z^q\left(t_0 + \frac{n}{q}\right)\right) d\langle W, W\rangle(s)$$
$$+ \int_{t_0+n/q}^t w\left(s, Z^q\left(t_0 + \frac{n}{q}\right)\right) dW(s). \quad (8)$$

Moreover, $Z^q(\cdot)$ can be resolved by stepwise iterations on the intervals $[t_0, t_0 + 1/q), [t_0 + 1/q, t_0 + 2/q), \ldots$ as follows. For $t \in [t_0, t_0 + 1/q]$

$$Z^q(t) = Z^q(t_0) + \int_{t_0}^t g(s, Z^q(t_0))\,ds$$
$$+ \int_{t_0}^t h(s, Z^q(t_0))\,d\langle W, W\rangle(s)$$
$$+ \int_{t_0}^t w(s, Z^q(t_0))\,dW(s),$$

and for $t \in [t_0 + 1/q, t_0 + 2/q]$

$$Z^q(t) = Z^q\left(t_0 + \frac{1}{q}\right)$$
$$+ \int_{t_0+1/q}^t g\left(s, Z^q\left(t_0 + \frac{1}{q}\right)\right) ds$$
$$+ \int_{t_0+1/q}^t h\left(s, Z^q\left(t_0 + \frac{1}{q}\right)\right) d\langle W, W\rangle(s)$$
$$+ \int_{t_0+1/q}^t w\left(s, Z^q\left(t_0 + \frac{1}{q}\right)\right) dW(s),$$

and so on. Defining $Z^q(t_0) = \tilde{Z}^q(t_0) = \tilde{Z}_0$ and $Z^q(t_0 + n/q) = \tilde{Z}^q(t)$ for $t_0 + n/q < t \leq (t_0 + (n+1)/q) \wedge T$, $n = 0, 1, 2, \ldots$, Equation (8) takes the following form

$$Z^q(t) = \tilde{Z}_0 + \int_{t_0}^t g(s, \tilde{Z}^q(s))\,ds$$
$$+ \int_{t_0}^t h(s, \tilde{Z}^q(s))\,d\langle W, W\rangle(s)$$
$$+ \int_{t_0}^t w(s, \tilde{Z}^q(s))\,dW(s), \quad (9)$$

for $t_0 \leq t \leq T$. Next, we derive an important result, which shows that for each $q \geq 1$, $\{Z^q(t)\}_{t \in [t_0, T]}$ is a well defined sequence in $M_G^2([t_0, T]; \mathbb{R}^n)$.

**Lemma 3.1:** *Let conditions (3) and (6) holds. For every $q \geq 1$ and any $T > 0$,*

$$\sup_{t_0 \leq t \leq T} \mathbb{E}(|Z^q(t)|^2) \leq K, \quad (10)$$

*where $K = G_1 e^{G_2 T}$, $G_1 = 4E|Z_0|^2 + 16T(T+2)(M+C)$, $G_2 = 16D(T+2)$ and $M, C, D$ are already defined positive constants.*

**Proof:** In view of the inequality $|\sum_{i=1}^4 a_i|^2 \leq 4 \sum_{i=1}^4 |a_i|^2$, Equation (9) gives

$$|Z^q(t)|^2 \leq 4|Z_0|^2 + 4\left|\int_{t_0}^t g(s, \tilde{Z}^q(s))\,ds\right|^2$$
$$+ 4\left|\int_{t_0}^t h(s, \tilde{Z}^q(s))\,d\langle W, W\rangle(s)\right|^2$$
$$+ 4\left|\int_{t_0}^t w(s, \tilde{Z}^q(s))\,dW(s)\right|^2$$

Apply subexpectation on both sides. Then in view of the Holder inequality (2.4), Lemmas 2.6 and 2.7, we proceed

as the following

$$E\left[\sup_{0\leq s\leq t}|Z^q(s)|^2\right]$$
$$\leq 4E|Z_0|^2 + 4T\int_{t_0}^t E|g(s,\tilde{Z}^q(s))|^2\,ds$$
$$+ 4T\int_{t_0}^t E|h(s,\tilde{Z}^q(s))|^2\,ds$$
$$+ 16\int_{t_0}^t E|w(s,\tilde{Z}^q(s))|^2\,ds$$
$$\leq 4E|Z_0|^2 + 8T\int_{t_0}^t E[|g(s,\tilde{Z}^q(s))$$
$$- g(s,0)|^2 + |g(s,0)|^2]\,ds$$
$$+ 8T\int_{t_0}^t E[|h(s,\tilde{Z}^q(s)) - h(s,0)|^2 + |h(s,0)|^2]\,ds$$
$$+ 32\int_{t_0}^t E[|w(s,\tilde{Z}^q(s)) - w(s,0)|^2 + |w(s,0)|^2]\,ds.$$

Utilizing conditions (3) and (6), the last inequality yields

$$E\left[\sup_{0\leq s\leq t}|Z^q(s)|^2\right]$$
$$\leq 4E|Z_0|^2 + 8T^2M + 8T^2M + 32TM$$
$$+ 8T\int_{t_0}^t E[\Psi(|\tilde{Z}^q(s)|^2)]\,ds$$
$$+ 8T\int_{t_0}^t E[\Psi(|\tilde{Z}^q(s)|^2)]\,ds$$
$$+ 32\int_{t_0}^t E[\Psi(|\tilde{Z}^q(s)|^2)]\,ds$$
$$\leq 4E|Z_0|^2 + 8T^2M + 8T^2M + 32TM$$
$$+ 16(T+2)\int_{t_0}^t E[\Psi(|\tilde{Z}^q(s)|^2)]\,ds$$
$$\leq 4E|Z_0|^2 + 8T^2M + 8T^2M + 32TM$$
$$+ 16(T+2)\int_{t_0}^t E[C + D|\tilde{Z}^q(s)|^2]\,ds$$
$$\leq 4E|Z_0|^2 + 8T^2M + 8T^2M + 32TM$$
$$+ 16CT(T+2) + 16D(T+2)\int_{t_0}^t E[|\tilde{Z}^q(s)|^2]\,ds.$$

In view of the notion $\tilde{Z}^q(s)$, we obtain

$$E\left[\sup_{t_0\leq s\leq t}|Z^q(s)|^2\right]$$
$$\leq 4E|Z_0|^2 + 8T^2M + 8T^2M + 32TM + 16CT(T+2)$$
$$+ 16D(T+2)\int_{t_0}^t E\left[\sup_{t_0\leq r\leq s}|Z^q(r)|^2\right]\,ds$$
$$\leq G_1 + G_2\int_{t_0}^t E\left[\sup_{t_0\leq r\leq s}|Z^q(r)|^2\right]\,ds,$$

where $G_1 = 4E|Z_0|^2 + 16T(T+2)(M+C)$ and $G_2 = 16D(T+2)$. Consequently, an application of Gronwall's inequality gives

$$E\left[\sup_{t_0\leq s\leq t}|Z^q(s)|^2\right] \leq G_1\,e^{G_2(t-t_0)},$$

which by assuming $t = T$ provides

$$E\left[\sup_{t_0\leq s\leq T}|Z^q(s)|^2\right] \leq G_1\,e^{G_2(T-t_0)} \leq G_1\,e^{G_2 T} = K.$$

The proof is complete. ∎

**Remark 3.2:** Lemma 3.1 shows that for every $q \geq 1$, $Z^q(t)$ is bounded in $M_G^2([t_0, T]; \mathbb{R}^d)$. By an identical way as used in lemma 3.1, one can prove that for any $T > 0$,

$$\sup_{t_0\leq t\leq T} E[|Z(t)|^2] \leq K, \quad (11)$$

where $K$ is a positive constant.

## 4. Estimates for the difference between approximate and exact solutions to SDEs in the G-framework

We now derive an important lemma, which will be utilized in the next theorem. Here we present estimates for the difference between an exact and approximate solutions for SDEs in the G-framework.

**Lemma 4.1:** *Assume that the hypothesis of Lemma 3.1 hold. Let $t_0 \leq r < t \leq T$. For all $q \geq 1$,*

$$\tilde{\mathbb{E}}[|Z^q(t) - Z^q(r)|^2] \leq H_1(t-r), \quad (12)$$

*where $H_1 = 12(T+2)(M + C + KD)$ and M,C,D,K are already defined positive constants.*

**Proof:** Let $t_0 \leq r < t \leq T$. For any $q \geq 1$, Equation (9) becomes

$$Z^q(t) - Z^q(r) = \int_r^t g(s,\tilde{Z}^q(s))\,ds$$
$$+ \int_r^t h(s,\tilde{Z}^q(s))\,d\langle W,W\rangle(s)$$
$$+ \int_r^t w(s,\tilde{Z}^q(s))\,dW(s).$$

Use the inequality $|\sum_{i=1}^3 a_i|^2 \leq 4\sum_{i=1}^3 |a_i|^2$ and apply subexpectation on both sides. Then in view of the Holder inequality (2.4), Lemmas 2.6 and 2.7, we proceed

as the following

$$E\left[\sup_{r\leq v<u\leq t}|Z^q(u)-Z^q(v)|^2\right]$$

$$\leq 3E\left[\sup_{r\leq v<u\leq t}\left|\int_v^u g(s,\tilde{Z}^q(s))\,ds\right|^2\right]$$

$$+3E\left[\sup_{r\leq v<u\leq t}\left|\int_v^u h(s,\tilde{Z}^q(s))\,d\langle W,W\rangle(s)\right|^2\right]$$

$$+3\tilde{\mathbb{E}}\left[\sup_{r\leq v<u\leq t}\left|\int_v^u w(s,\tilde{Z}^q(s))\,dW(s)\right|^2\right]$$

$$\leq 3T\int_r^t E[|g(s,\tilde{Z}^q(s))|^2]\,ds$$

$$+3T\int_r^t E[|h(s,\tilde{Z}^q(s))|^2]\,ds$$

$$+12\int_r^t E[|w(s,\tilde{Z}^q(s))|^2]\,ds$$

$$\leq 6T\int_r^t E[|g(s,\tilde{Z}^q(s))-g(s,0)|^2+|g(s,0)|^2]\,ds$$

$$+6T\int_r^t E[|h(s,\tilde{Z}^q(s))-h(s,0)|^2+|h(s,0)|^2]\,ds$$

$$+24\int_r^t E[|w(s,\tilde{Z}^q(s))-w(s,0)|^2+|w(s,0)|^2]\,ds$$

$$\leq 6MT(t-r)+6MT(t-r)+24M(t-r)$$

$$+12(T+2)\int_r^t E[\Psi(|\tilde{Z}^q(s)|^2)]\,ds$$

$$\leq 6MT(t-r)+6MT(t-r)+24M(t-r)$$

$$+12(T+2)\int_r^t E[C+D|\tilde{Z}^q(s)|^2]\,ds$$

$$\leq 6MT(t-r)+6MT(t-r)+24M(t-r)$$

$$+12C(T+2)(t-r)$$

$$+12D(T+2)\int_r^t E[|\tilde{Z}^q(s)|^2]\,ds$$

In view of the notion of $\tilde{Z}^q(s)$, we have

$$E\left[\sup_{r\leq v<u\leq t}|Z^q(u)-Z^q(v)|^2\right]$$

$$\leq 6MT(t-r)+6MT(t-r)$$

$$+24M(t-r)+12C(T+2)(t-r)$$

$$+12D(T+2)\int_r^t E\left[\sup_{t_0\leq v\leq s}|Z^q(v)|^2\right]ds$$

By virtue of Lemma 3.1, we get

$$E\left[\sup_{r\leq v<u\leq t}|Z^q(u)-Z^q(v)|^2\right]$$

$$\leq 6MT(t-r)+6MT(t-r)+24M(t-r)$$

$$+12C(T+2)(t-r)$$

$$+12D(T+2)K(t-r).$$

Consequently,

$$E[|Z^q(t)-Z^q(r)|^2]\leq H_1(t-r),$$

where $H_1=12(T+2)(M+C+KD)$. The proof stands completed. ■

**Remark 4.2:** Using identical arguments as used in Lemma 4.1, one can prove that

$$E[|X(t)-X(r)|^2]\leq H_1(t-r), \qquad (13)$$

where $H_1$ is a positive constant.

**Theorem 4.3:** Let (3) and (6) hold. Then for all $q\geq 1$ and any $T>0$,

$$E\left[\sup_{t_0\leq s\leq T}|Z^q(s)-Z(s)|^2\right]$$
$$\leq 6T(T+2)\left[C+\frac{2DH_1}{q}\right]e^{12(T+2)D(t-t_0)},$$

where $C$, $D$ and $H_1$ are positive constants.

**Proof:** Using the fundamental inequality $|\sum_{i=1}^3 a_i|^2 \leq 4\sum_{i=1}^3 |a_i|^2$, from (2) and (9) we derive

$$|Z^q(t)-Z(t)|^2$$

$$\leq 3\left|\int_{t_0}^t [g(s,\tilde{Z}^q(s))-g(s,Z(s))]\,ds\right|^2$$

$$+3\left|\int_{t_0}^t [h(s,\tilde{Z}^q(s))-h(s,Z(s))]\,d\langle W,W\rangle(s)\right|^2$$

$$+3\left|\int_{t_0}^t [w(s,\tilde{Z}^q(s))-w(s,Z(s))]\,dW(s)\right|^2.$$

Apply subexpectation on both sides. Then in virtue of the Holder inequality (2.4), Lemmas 2.6 and 2.7, we derive

$$E\left[\sup_{t_0\leq s\leq t}|Z^q(s)-Z(s)|^2\right]$$

$$\leq 3T\int_{t_0}^t E[|g(s,\tilde{Z}^q(s))-g(s,Z(s))|^2]\,ds$$

$$+3T\int_{t_0}^t E[|h(s,\tilde{Z}^q(s))-h(s,Z(s))|^2]\,ds$$

$$+12\int_{t_0}^t E[|w(s,\tilde{Z}^q(s))-w(s,Z(s))|^2]\,ds.$$

Applying the non-uniform Lipschitz condition we have

$$E\left[\sup_{t_0 \leq s \leq t} |Z^q(s) - Z(s)|^2\right]$$
$$\leq 6(T+2) \int_{t_0}^{t} E[\Upsilon(|\tilde{Z}^q(s) - Z(s)|^2)] \, ds$$
$$\leq 6T(T+2)C + 6(T+2)D \int_{t_0}^{t} E[|\tilde{Z}^q(s) - Z(s)|^2] \, ds$$
$$= 6T(T+2)C + 6(T+2)D$$
$$\times \int_{t_0}^{t} E[|\tilde{Z}^q(s) - \tilde{Z}(s) + \tilde{Z}(s) - Z(s)|^2] \, ds$$
$$\leq 6T(T+2)C + 12(T+2)D \int_{t_0}^{t} E[|\tilde{Z}^q(s) - \tilde{Z}(s)|^2] \, ds$$
$$+ 12(T+2)D \int_{t_0}^{t} E[|\tilde{Z}(s) - Z(s)|^2] \, ds$$
$$= 6T(T+2)C + 12(T+2)D$$
$$\times \int_{t_0}^{t} E[|\tilde{Z}^q(s) - \tilde{Z}(s)|^2] \, ds + \mathcal{N},$$

where $\mathcal{N} = 12(T+2)D \int_{t_0}^{t} E[|\tilde{Z}(s) - Z(s)|^2] \, ds$. By an application of the Grownwall's inequality we derive

$$E\left[\sup_{t_0 \leq s \leq t} |Z^q(s) - Z(s)|^2\right]$$
$$\leq [6T(T+2)C + \mathcal{N}] e^{12(T+2)D(t-t_0)}. \quad (14)$$

Using Lemma 4.1, we estimate $\mathcal{N}$ as follows

$$\mathcal{N} = 12(T+2)D \sum_{n \geq 0} \int_{t_0 + n/q}^{(t_0 + (n+1)/q) \wedge T}$$
$$\times E\left[\left|Z\left(t_0 + \frac{n}{q}\right) - Z(s)\right|^2\right] ds$$
$$\leq 12T(T+2)DH_1 \frac{1}{q},$$

substituting the value of $\mathcal{N}$ in (14) provides,

$$E\left[\sup_{t_0 \leq s \leq t} |Z^q(s) - Z(s)|^2\right]$$
$$\leq \left[6T(T+2)C + 12T(T+2)DH_1 \frac{1}{q}\right] e^{12(T+2)D(t-t_0)}$$
$$= 6T(T+2) \left[C + \frac{2DH_1}{q}\right] e^{12(T+2)D(t-t_0)}.$$

Consequently, by assuming $t = T$,

$$E\left[\sup_{t_0 \leq s \leq T} |Z^q(s) - Z(s)|^2\right]$$
$$\leq 6T(T+2) \left[C + \frac{2DH_1}{q}\right] e^{12(T+2)D(T-t_0)}.$$

The proof stands completed. ∎

## 5. Conclusion

In recent years, the importance of SDEs has become more apparent due to their applications in modelling real life phenomena. Subject to the Lipschitz conditions, the existence theory for stochastic functional differential equations (SFDEs) in the G-framework was developed by Ren, Bi and Sakthivel [19]. The stated theory was extended to non-uniform Lipschitz conditions by Faizullah [20–22] and to discontinuous coefficients by Faizullah, Rahman, Afzal and Chohan [23]. Further, Faizullah established the pth moment estimates for the stated equations [24,25]. It is expected that the techniques used in the present paper can be used in several different directions such as to find estimates for the difference between the exact and approximate solutions for the above stated SFDEs in the G-framework, neural stochastic differential equations driven by G-Brownian motion [26] and stochastic differential equations with piecewise arguments in the G-framework etc. We hope that the current study will play a key role to establish a framework for the above mentioned problems.

## Disclosure statement

No potential conflict of interest was reported by the authors.

## Funding

This work is supported by the Commonwealth Scholarship Commission in the United Kingdom with CSC ID: PKRF-2017-429.

## ORCID

*Faiz Faizullah*  http://orcid.org/0000-0001-7474-7754
*Ilyas Khan*  http://orcid.org/0000-0002-2056-9371
*Mukhtar M. Salah*  http://orcid.org/0000-0002-2736-8134
*Ziyad Ali Alhussain*  http://orcid.org/0000-0001-8593-0239

## References

[1] Zamana G, Kang YH, Jung IH. Stability analysis and optimal vaccination of an SIR epidemic model. BioSystems. 2008;93:240–249.
[2] Zaman G, Kang YH, Cho G, et al. Optimal strategy of vaccination & treatment in an SIR epidemic model. Math Comput Simul. 2017;136:63–77.
[3] Khan L, Ali F, Shah NA. Interaction of magnetic field with heat and mass transfer in free convection flow of a Walters'-B fluid. Eur Phys J Plus. 2016;131(77):1–15.
[4] Peng S. Multi-dimentional G-Brownian motion and related stochastic calculus under G-expectation. Stoch Process Appl. 2008;12:2223–2253.
[5] Peng S. G-expectation, G-Brownian motion and related stochastic calculus of Ito's type. In: Benth et al., editors. The Abel symposium 2005, Abel symposia 2. Berlin: Springer-verlag; 2006. p. 541–567.
[6] Gao F. Pathwise properties and homeomorphic flows for stochastic differential equations driven by G-Brownian motion. Stoch Process Appl. 2009;10:3356–3382.
[7] Faizullah F. A note on the Caratheodory approximation scheme for stochastic differential equations under

G-Brownian motion. Z Naturforschung A. 2012;67a: 699–704.
[8] Bai X, Lin Y. On the existence and uniqueness of solutions to stochastic differential equations driven by G-Brownian motion with Integral-Lipschitz coefficients. Acta Math Appl Sin. 2014;30(3):589–610.
[9] Faizullah F. Existence of solutions for stochastic differential equations under G-Brownian motion with discontinuous coefficients. Z Naturforschung A. 2012;67a: 692–698.
[10] Denis L, Hu M, Peng S. Function spaces and capacity related to a sublinear expectation: application to G-Brownian motion paths. Potential Anal. 2010;34: 139–161.
[11] Ullah R, Faizullah F. On existence and approximate solutions for stochastic differential equations in the framework of G-Brownian motion. Eur Phys J Plus. 2017;132:435–443.
[12] Gao F, Jiang H. Approximation Theorem for Stochastic Differential Equations Driven by G-Brownian Motion. Prog Probab. 2011;65:73–81.
[13] Gu Y, Ren Y, Sakthivel R. Square-mean pseudo almost automorphic mild solutions for stochastic evolution equations driven by G-Brownian motion. Stoch Anal Appl. 2016;34(3):528–545.
[14] Li X, Peng S. Stopping times and related Ito's calculus with G-Brownian motion. Stoch Process Appl. 2011;121:1492–1508.
[15] Luo P, Wang F. Stochastic differential equations driven by G-Brownian motion and ordinary differential equations. Stoch Process Appl. 2014;124:3869–3885.
[16] Song Y. Properties of hitting times for G-martingale and their applications. Stoch Process Appl. 2011;121(8): 1770–1784.
[17] Xua J, Zhang B. Martingale characterization of G-Brownian motion. Stoch Process Appl. 2009;119: 232–248.
[18] Mao X. Stochastic differentail equations and their applications. Chichester, England: Horwood Publishing; 1997.
[19] Ren Y, Bi Q, Sakthivel R. Stochastic functional differential equations with infinite delay driven by G-Brownian motion. Math Methods Appl Sci. 2013;36(13):1746–1759.
[20] Faizullah F. Existence and uniqueness of solutions to SFDEs driven by G-Brownian motion with non-Lipschitz conditions. J Comput Anal Appl. 2017;2(23):344–354.
[21] Faizullah F. Stochastic functional differential equations driven by G-Browniain motion with monotone nonlinearity. 2018. Available from: https://arxiv.org/abs/1806.07824.
[22] Faizullah F. Existence and asymptotic properties for the solutions to nonlinear SFDEs driven by G-Brownian motion with infinite delay. 2018. Available from: https://arxiv.org/abs/1805.10658.
[23] Faizullah F, Rehman M-U, Shehzad M, et al. On existence and comparison results for solutions to stochastic functional differential equations in the G-framework. J Comput Anal Appl. 2017;4(23):693–702.
[24] Faizullah F. On the pth moment estimates of solutions to stochastic functional differential equations in the G-framework. SpringerPlus. 2016;5(872):1–11.
[25] Faizullah F. A note on pth moment estimates for stochastic functional differential equations in the framework of G-Brownian motion,Iranian Journal of Science and Technology. Trans. A. 2016;3(40):1–8.
[26] Faizullah F. Existence results and moment estimates for NSFDEs driven by G-Brownian motion. J Comput Theor Nanosci. 2016;7(13):1–8.